\documentclass[a4paper,11pt]{amsart}
\usepackage{amssymb,amsfonts,amsxtra,amscd}
\usepackage[all]{xy}
\usepackage{fullpage}
\theoremstyle{plain}
\newtheorem{theorem}{Theorem}[section]
\newtheorem{lemma}[theorem]{Lemma}
\newtheorem{cor}[theorem]{Corollary}
\newtheorem{prop}[theorem]{Proposition}
\theoremstyle{definition}
\newtheorem{defi}[theorem]{Definition}
\theoremstyle{remark}
\newtheorem{rem}[theorem]{Remark}
\numberwithin{equation}{section}
\newcommand{\Ai}{$A_\infty$}
\DeclareMathOperator{\Hom}{Hom}
\DeclareMathOperator{\Obs}{Obs}
\DeclareMathOperator{\End}{End}
\begin{document}
\begin{abstract}
In this paper we will study deformations of \Ai-algebras. We will also answer questions relating to Moore algebras which are one of the simplest nontrivial examples of an \Ai-algebra. We will compute the Hochschild cohomology of odd Moore algebras and classify them up to a unital weak equivalence. We will construct miniversal deformations of particular Moore algebras and relate them to the universal odd and even Moore algebras. Finally we will conclude with an investigation of formal one-parameter deformations of an \Ai-algebra.
\end{abstract}
\title{On the Classification of Moore Algebras and their Deformations}
\author{Alastair Hamilton}
\address{Mathematics Department, Bristol University, Bristol, England. BS8 1TW.}
\email{a.e.hamilton@bristol.ac.uk}
%\date{$19^{\mathbf{th}}$ February 2003}
\thanks{The work of the author was supported by an EPSRC grant No. GR/SO7148/01}
\maketitle
\section{Introduction}
The notion of an \Ai-algebra, first introduced in \cite{Stasheff} by Stasheff, is much the same as that of an ordinary DGA (differential graded algebra), except that it possesses extra structure in the form of higher multiplications. It is this extra structure however that proves to be useful when it comes to making explicit computations. In section \ref{Preliminaries} I shall briefly introduce the notion of an \Ai-algebra and define its Hochschild cohomology. Further details can be found in \cite{Andrey} and \cite{Bernhard}.

The first examples of Moore algebras were seen in Kontsevich's paper \cite{Kont} where they were related to Morita-Miller-Mumford classes in the cohomology of moduli spaces of algebraic curves. Moore algebras were introduced in their full generality in \cite{Andrey} by Lazarev, where some results for even Moore algebras were presented. In this paper I shall present some results for odd Moore algebras and deformations of Moore algebras in general. Moore algebras will be introduced in section \ref{Preliminaries} but the reader should refer to \cite{Andrey} for a more detailed description.

Throughout the rest of this paper we will be working over an evenly graded and commutative ring $R$. A Moore algebra of degree $d$ is a (unital) \Ai-structure on the $R$-module,
\begin{displaymath} A:=\Sigma^{d+1}R \oplus R \end{displaymath}
which consists of a generator $1$ in degree $0$ and a generator $y$ in degree $d+1$. This has an underlying two cell complex of the form $\Sigma^d R \to R$. A Moore algebra is then even (odd) if $d$ is even (odd). An \Ai-structure on $A$ is determined by a coderivation $m$ on the tensor coalgebra $T\Sigma A$ such that $m^2=0$ (cf. \cite{Andrey}, \cite{Bernhard}). Remarkably, in the case of Moore algebras, the condition $m^2=0$ places no restrictions on our choice of coderivation (cf. \cite{Andrey}).

In section \ref{Classification} I shall present the classification result for odd Moore algebras. The result is similar to that obtained in the even case in \cite{Andrey}, although there is no initially apparent reason as to why this should be so. In section \ref{Cohomology} I shall calculate the Hochschild cohomology of odd Moore algebras, the even case having been treated already in \cite{Andrey}. The problem of calculating the Hochschild cohomology of odd Moore algebras was posed in \cite{Andrey} by Lazarev. The classification result in section \ref{Classification} allows us to transform the odd Moore algebra to one in which the Hochschild cohomology is easier to compute.

Deformation theory for associative algebras was first developed in \cite{Murray} by Gerstenhaber and the deformation theory for \Ai-algebras can be seen as a direct generalisation of this. In section \ref{Deform} I shall describe the deformation theory for \Ai-algebras analogous to that of \cite{Fuchs}, \cite{Alicemini} and \cite{Penkava}. Essentially this involves deforming the \Ai-structure (a derivation on the cobar construction) according to the methods of \cite{dgadef}. In \cite{Penkava} the link between Hochschild cohomology and deformation theory was made and this link will provide the motivation for section \ref{Mini}.

In section \ref{Mini} I shall describe the ``miniversal'' deformations of the ``trivial'' Moore algebra $R[X]/X^2$ which is the Moore algebra that corresponds to the coderivation which sends $[y]^{\otimes i}$ to $0$. According to \cite{Schlessinger}, universal deformations correspond to the pro-representability of a certain functor whilst miniversal deformations correspond to the hull of that functor, however, the approach taken in this paper will be slightly different.

In section \ref{Oneparam} we will look at formal one-parameter deformations of \Ai-algebras and discuss some of the theory involved. The theory of formal one-parameter deformations of an \Ai-algebra is closely related to its Hochschild cohomology. This theory was first presented in the context of ordinary DGA's by Gerstenhaber and Wilkerson in \cite{dgadef}. The corresponding theory for \Ai-algebras was developed by Penkava and Weldon in \cite{Weldon} and Penkava and Fialowski in \cite{Alicemini} (although Penkava and Fialowski did not explicitly consider one-parameter deformations). Whilst this section contains some results not seen in \cite{Weldon} and \cite{Alicemini} which are analogues of results contained in \cite{dgadef}, it is hoped that by simply translating the material into the language of the cobar construction, the ideas involved will stand out more clearly.

\emph{Acknowledgement:} The author would like to thank Andrey Lazarev for his advice and assistance which was instrumental to the completion of this work. The author would also like to thank Murray Gerstenhaber for a helpful consultation.
\section{\Ai-Algebras and Hochschild Cohomology} \label{Preliminaries}
In this section I shall briefly introduce the notion of an \Ai-algebra and a Moore algebra and define its Hochschild cohomology.

\begin{defi} \label{def_AI}
A unital \Ai-algebra is a free graded module $A$ defined over an evenly graded and commutative ring $R$ which has an $R$-basis $\{1,y_i \ i \in I\}$. An \Ai-structure on $A$ is then a continuous derivation of degree $-1$ on the $R$-algebra of formal power series $R[[\tau,t_1,t_2,\ldots]]$ whose square is zero. The generators $\{\tau,\boldsymbol{t}\}=\tau,t_1,t_2,\ldots$ have degrees $|\tau|=-1$, $|t_i|=-|y_i|-1$. Define the map $ad\tau$ by $ad\tau(x):=[\tau,x]$, then a unital \Ai-structure is a derivation of the form,
\begin{equation} \label{eqn_AI}
m=A(\boldsymbol{t})\partial_\tau + \sum_{i \in I} B_i(\boldsymbol{t})\partial_{t_i} +ad\tau -\tau^2\partial_\tau
\end{equation}
such that $m^2=0$ and the power series $A(\boldsymbol{t})$, $B_i(\boldsymbol{t}),i \in I$ have vanishing constant terms. In this paper we will only be considering unital \Ai-algebras.
\end{defi}

\begin{defi} \label{def_AImorph}
Suppose we have two \Ai-algebras $A$ and $A'$ with \Ai-structures $m$ and $m'$ respectively. An \Ai-morphism of $A$ into $A'$ is a continuous morphism of $R$-algebras $f:R[[\tau',\boldsymbol{t}']] \to R[[\tau,\boldsymbol{t}]]$ of degree $0$ such that $mf=fm'$. Such a morphism is determined by its action on the generators and is therefore given by a collection of power series $(G(\tau,\boldsymbol{t}),F_1(\tau,\boldsymbol{t}),F_2(\tau,\boldsymbol{t}),\ldots)$:
\begin{displaymath}
\begin{array}{ccc}
\tau' & \mapsto & G(\tau,\boldsymbol{t}) \\
t_1' & \mapsto & F_1(\tau,\boldsymbol{t}) \\
\vdots & \mapsto & \vdots 
\end{array}
\end{displaymath}
A unital \Ai-isomorphism $f$ is specified by a collection of power series with vanishing constant terms $(\tau+G(\boldsymbol{t}),F_1(\boldsymbol{t}),F_2(\boldsymbol{t}),\ldots)$ such that $mf=fm'$. In addition the collection $\boldsymbol{F}(\boldsymbol{t})=(F_1(\boldsymbol{t}),F_2(\boldsymbol{t}),\ldots)$ determines a map $\boldsymbol{F}:R[[\boldsymbol{t}]] \to R[[\boldsymbol{t}]]$ which must be invertible. In this paper we will only be considering unital \Ai-isomorphisms.
\end{defi}

\begin{rem} \label{rem_hmult}
There is an alternative definition of an \Ai-algebra in terms of a collection of higher multiplications $m_i:A^{\otimes i} \to A, \ i\geq 1$. The condition $m^2=0$ above translates to imposing certain restrictions on the $m_i$'s, one of which is that the map $m_1$ is a graded derivation with respect to the multiplication $m_2$ and that $m_1$ gives $A$ the structure of a differential graded module.  An \Ai-morphism is similarly given by a collection of maps $f_i:A^{\otimes i} \to A', \ i\geq 1$ satisfying certain conditions, one of which is that $f_1:A \to A'$ is a map of differential graded modules with respect to the differential $m_1$. This morphism of \Ai-algebras is an isomorphism if and only if $f_1$ is an isomorphism. A morphism of \Ai-algebras is called a weak equivalence if $f_1$ induces an isomorphism in the homology of the differential graded modules $A$ and $A'$.

This alternative definition derives from the fact that $R[[\tau,\boldsymbol{t}]]$ is the $R$-linear dual to the tensor coalgebra $T\Sigma A$. Continuous derivations of $R[[\tau,\boldsymbol{t}]]$ are in one to one correspondence with coderivations of $T\Sigma A$ which in turn are specified by the maps $m_i:A^{\otimes i} \to A$ (cf. \cite{Andrey}).
\end{rem}

\begin{defi} \label {def_MA}
A Moore algebra is an \Ai-algebra with an underlying module structure $A:=\Sigma^{d+1}R \oplus R$. The algebra is even (odd) if $d$ is even (odd). This algebra has a basis $1,y$ where $|y|=d+1$ and $|1|=0$. An \Ai-structure is then a continuous derivation $m$ of $R[[\tau,t]]$, where $|\tau|=-1$ and $|t|=-(d+2)$. Remarkably the condition $m^2=0$ places no restrictions on the choice of derivation. This means that even Moore algebras are characterised by formal power series $u(t) \in R[[t]]$ of degree $-2$;
\begin{equation} \label{eqn_evendef}
m=\sum_{i=1}^\infty u_it^i\partial_\tau +ad\tau-\tau^2\partial_\tau
\end{equation}
whilst odd Moore algebras are characterised by pairs of formal power series $v(t),w(t) \in R[[t]]$ of degree $-(d+3)$ and $-2$ respectively;
\begin{equation} \label{eqn_odddef}
m=\sum_{i=1}^{\infty} v_it^{2i}\partial_t + \sum_{i=1}^{\infty} w_it^{2i}\partial_{\tau} + ad\tau -{\tau}^2\partial_{\tau}
\end{equation}
\end{defi}

\begin{rem}
A trivial but important observation is that the parameter $t$ is even (odd) for even (odd) Moore algebras. If the parameter $t$ is odd then a formal power series is even (odd) if and only if it consists entirely of even (odd) powers of $t$.
\end{rem}

We will also need a way to calculate the Hochschild cohomology of an \Ai-algebra. Hochschild cohomology is closely related to deformation theory for \Ai-algebras and an account of this is given in \cite{Penkava}. For this reason we will be chiefly concerned with the Hochschild cohomology of an \Ai-algebra with coefficients in itself. Further details on the Hochschild cohomology of \Ai-algebras can be found in \cite{Andrey}.

\begin{defi} \label{def_Hoch}
The Hochschild cohomology of an \Ai-algebra $A$ with coefficients in $A$ is defined via the complex $C^*(A,A)$ of normalised derivations of the algebra $R[[\tau,\boldsymbol{t}]]$. A normalised derivation is one of the form,
\[\xi=A(\boldsymbol{t})\partial_\tau + \sum_{i \in I}B_i(\boldsymbol{t})\partial_{t_i} \]
where $A(\boldsymbol{t})$, $B_i(\boldsymbol{t}), i \in I$ have vanishing constant terms. The \Ai-structure $m$ gives rise to a differential on this complex which is given by $d(\xi):=[\xi,m]$. The Hochschild cohomology of $A$ with coefficients in $A$ is then the cohomology of this complex $C^*(A,A)$.
\end{defi}

\begin{rem}
The complex $C^*(A,A)$ is graded according to the degrees of the maps $\xi:R[[\tau,\boldsymbol{t}]] \to R[[\tau,\boldsymbol{t}]]$. I will refer to this as the standard grading and it is this grading that must be used for the purpose of calculations. It will be assumed throughout the rest of the paper that we are using the standard grading unless an explicit reference is made to the contrary.

There is however an alternative grading that is more customarily used when discussing deformation theory. Continuous derivations of the algebra $R[[\tau,\boldsymbol{t}]]$ are in one to one correspondence with coderivations of the tensor coalgebra $T \Sigma A$. These coderivations have a bi-grading,
\[ C^i_j=\{f \in \Hom_R(A^{\otimes i},A): |f|=j \} \]
The associated total grading $T_k:=\prod_{i-j=k}C^i_j$ forms an alternative grading for $C^*(A,A)$ which I shall refer to as the classical grading. If a map $\xi$ has degree $d$ in the standard grading then its corresponding degree in the classical grading is $1-d$. One of the consequences of this is that infinitesimals live in even components of the classical grading, so that the deformation theory of \Ai-algebras is consistent with classical deformation theory.
\end{rem}
\section{Classification of Odd Moore Algebras} \label{Classification}
As mentioned earlier an odd Moore algebra is characterised by a pair of power series $v(t),w(t) \in R[[t]]$ and has the \Ai-structure given by \eqref{eqn_odddef}. We would like to classify these structures up to a weak equivalence. Since the Moore algebra is odd, the multiplication map $m_1=0$ and so weak equivalences are isomorphisms of \Ai-algebras. Such \Ai-isomorphisms are continuous automorphisms of the $R$-algebra $R[[\tau,t]]$ and are specified by a pair of odd power series with vanishing constant terms $G(t),F(t)\in R[[t]]$ (cf. definition \ref{def_AImorph}) where the coefficient of $t$ in $F(t)$ is invertible. The degrees of $G(t)$ and $F(t)$ are $-1$ and $-(d+2)$ respectively. Suppose we are given two pairs of power series $(G(t),F(t))$ and $(G'(t),F'(t))$ corresponding to continuous automorphisms of the algebra $R[[\tau,t]]$. The pair $(G(t)+G'(F(t)),F'(F(t)))$ corresponds to the composition $(G,F)\circ (G',F')$. Since the identity automorphism corresponds to the pair $(0,t)$, we conclude that given a pair of power series with vanishing constant terms $(G(t),F(t))$ corresponding to a continuous automorphism of the algebra $R[[\tau,t]]$, the inverse of this automorphism corresponds to the pair of power series $(-G(F^{-1}(t)),F^{-1}(t))$, where $F^{-1}$ is the inverse of the map $t \mapsto F(t)$. These continuous automorphisms of the algebra $R[[\tau,t]]$ act on the space of unital \Ai-structures by conjugation so that equivalence classes of unital \Ai-algebras correspond to orbits of this action.

Let me now describe a group $H \subset R[[t]]$ which consists of power series $F(t)$ which have degree $|t|=-(d+2)$ and a vanishing constant term. Furthermore the coefficient of $t$ in $F(t)$ must be an invertible element of $R$. The multiplication in $H$ is given by $(F(t),G(t)) \mapsto F(G(t))$. $H$ then acts on the right of the set consisting of formal power series in $R[[t]]$ which have degree -2 and a vanishing constant term by the formula $(u(t),f(t)) \mapsto u(f(t))$.

We shall now describe the unital weak equivalence classes of an odd Moore algebra in terms of this group action, however we will have to assume that the element $2 \in R$ is invertible.

\begin{theorem} \label{oddclass}
The set consisting of unital weak equivalence classes of odd Moore algebras of degree $d$ defined over the ring $R$ ($\tfrac{1}{2} \in R$) is in one to one correspondence with the set of orbits of the group $H$ acting on the set consisting of formal power series of degree $-2$ with a vanishing constant term, as described above. 
\end{theorem}

\begin{proof}
Let us determine how the continuous automorphisms of the $R$-algebra $R[[\tau,t]]$ as described above act on the elements $t^{2i}\partial_{\tau}$, $t^{2i}\partial_t$ and $m_0:=ad\tau -{\tau}^2\partial_{\tau}$.

\begin{equation} \label{eqn_deetau}
\begin{split}
& \begin{array}{lcl}
(G,F)\circ(t^{2i}\partial_{\tau})\circ(G,F)^{-1}(t) & = & 0 ; \\
(G,F)\circ(t^{2i}\partial_{\tau})\circ(G,F)^{-1}(\tau) & = & (G,F)\circ(t^{2i}\partial_{\tau})(\tau-G(F^{-1}(t))) \\
& = & (G,F)(t^{2i})=[F(t)]^{2i} ; \\
\end{array}
\\
& \text{therefore, } (G,F)\circ(t^{2i}\partial_{\tau})\circ(G,F)^{-1} = [F(t)]^{2i}\partial_{\tau}.
\end{split}
\end{equation}

\begin{equation} \label{eqn_deetee}
\begin{split}
& \begin{array}{lclr}
(G,F)\circ(t^{2i}\partial_t)\circ(G,F)^{-1}(t) & = & (G,F)\circ(t^{2i}\partial_t)[F^{-1}(t)] \\
 & = & (G,F)(t^{2i-1}F^{-1}(t)) & \text{ as $F(t)$ is odd} \\
 & = & t[F(t)]^{2i-1} ; \\
(G,F)\circ(t^{2i}\partial_t)\circ(G,F)^{-1}(\tau) & = & (G,F)\circ(t^{2i}\partial_t)[\tau -G(F^{-1}(t))] \\
& = & -(G,F)(t^{2i-1}G(F^{-1}(t))) & \text{ as $G(t)$ is odd} \\
& = & -G(t)[F(t)]^{2i-1} ; \\
\end{array}
\\
& \text{therefore, } (G,F)\circ(t^{2i}\partial_t)\circ(G,F)^{-1} = [F(t)]^{2i-1}(t\partial_t-G(t)\partial_{\tau}).
\end{split}
\end{equation}

\begin{equation} \label{eqn_deemin}
\begin{split}
& \begin{array}{lcl} 
(G,F)\circ m_0 \circ(G,F)^{-1}(t) & = & (G,F)\circ m_0[F^{-1}(t)] \\
& = & (G,F)([\tau,F^{-1}(t)]) \\
& = & [\tau,t] + 2tG(t) ; \\
(G,F)\circ m_0 \circ(G,F)^{-1}(\tau) & = & (G,F)\circ m_0[\tau -G(F^{-1}(t))] \\
& = & (G,F)(\tau^2-[\tau,G(F^{-1}(t))]) \\
& = & (\tau+G(t))^2-[\tau+G(t),G(t)] \\
& = & \tau^2-G(t)^2 ; \\
\end{array}
\\
& \begin{array}{lcl} 
\text{therefore, } (G,F)\circ m_0 \circ(G,F)^{-1} & = & ([\tau,t]+2tG(t))\partial_t + (\tau^2-G(t)^2)\partial_{\tau} \\
& = &  m_0 + 2tG(t)\partial_t - G(t)^2\partial_{\tau}.
\end{array}
\end{split}
\end{equation}

From the above equations we conclude that for evenly graded power series $v(t),w(t) \in R[[t]]$,
\begin{displaymath}
\begin{array}{lcl}
(G,F)\circ(w(t)\partial_{\tau})\circ(G,F)^{-1} & = & w(F(t))\partial_{\tau} \\
(G,F)\circ(v(t)\partial_t)\circ(G,F)^{-1} & = & \dfrac{v(F(t))}{F(t)}(t\partial_t-G(t)\partial_{\tau}) \\
\end{array}
\end{displaymath}
and obtain a formula for the action on \eqref{eqn_odddef}:
\begin{equation} \label{eqn_oddact}
\begin{split}
(G,F)\circ (m_0+v(t)\partial_t+w(t)\partial_{\tau})\circ(G,F)^{-1} = & m_0+\left\{2tG(t)+\dfrac{tv(F(t))}{F(t)}\right\}\partial_t \\
 & +\left\{-\dfrac{G(t)v(F(t))}{F(t)}+w(F(t))-G(t)^2\right\}\partial_{\tau} \\
\end{split}
\end{equation}

Now if we choose $F(t)=t$ and $G(t)=-\tfrac{v(t)}{2t}$ then \eqref{eqn_odddef} is sent to $m_0+\{\left(\scriptstyle\frac{v(t)}{2t}\right)^2+w(t)\}\partial_\tau$. We see that by using $G(t)$, every odd Moore algebra can be transformed to one of the form $m_0+u(t)\partial_\tau$. $F(t)$ then acts very simply on these forms sending $m_0+u(t)\partial_\tau$ to $m_0+u(F(t))\partial_\tau$. Equation \eqref{eqn_oddact} shows that two Moore algebras of the form $m_0+u(t)\partial_\tau$ can only be equivalent if $G(t)=0$, therefore weak equivalence classes correspond to orbits of the group $H$ acting on power series of degree $-2$ with a vanishing constant term.
\end{proof}

\begin{rem}
Using the change of variable $t \mapsto t^2$ we can see that this orbit space is the same as the orbit space of $H$ acting on power series with a vanishing constant term whose coefficient at $t^i$ has degree $2i(d+2)-2$. These power series consist of odd powers of $t$ as well as even.
\end{rem}

\begin{rem}
The orbit space of $H$ has already been described for fields and discrete valuation rings. The reader should refer to \cite[\S 6]{Andrey} for further details.
\end{rem}
\section{The Hochschild Cohomology of Odd Moore Algebras} \label{Cohomology}

We shall now turn our attention to a problem posed in \cite{Andrey}, that of calculating the Hochschild cohomology of an odd Moore algebra. As mentioned in section \ref{Preliminaries} the Hochschild cohomology of an odd Moore algebra is defined via the complex of normalised derivations of $R[[\tau,t]]$, the differential being given by $d(\xi):=[\xi,m]$ where $m$ is the \Ai-structure. These normalised derivations are of the form $\xi=A(t)\partial_\tau+B(t)\partial_t$ where $A(t)$ and $B(t)$ have vanishing constant terms. We will refer to the odd and even parts of $A(t)$ as $A_1(t)$ and $A_2(t)$ respectively. The \Ai-structure $m$ is specified by characteristic power series $v(t),w(t) \in R[[t]]$ as in \eqref{eqn_odddef}.

\begin{rem}
Concerning notation, whenever I write $a'(t)$ for a derivative, I am referring to the following:
\begin{displaymath}
a'(t):=\sum_{i=1}^\infty ia_it^{i-1}
\end{displaymath}
\end{rem}

\begin{lemma} \label{Hochformula}
Suppose we have a normalised derivation $\xi=A(t)\partial_\tau+B(t)\partial_t$, then we have the following formula for the differential $d$:
\begin{align*}
[\xi,m]= & \{B_1(t)w'(t)-A_1(t)\tfrac{v(t)}{t}\}\partial_\tau \\
 & +B_1(t)\{v'(t)-\tfrac{v(t)}{t}\}\partial_t \\
 & +2tA_1(t)\partial_t
\end{align*}
\end{lemma}

\begin{proof}
The contribution from $w(t)\partial_\tau$ is given by,
\begin{equation}
\begin{split}
& \begin{array}{lcl}
[\xi,w(t)\partial_\tau](t) & = & 0 ; \\
\\
\xi \: w(t)\partial_\tau(\tau) & = & \xi(w(t))=B(t)\partial_t(w(t))=B_1(t)w'(t)\\
w(t)\partial_\tau \: \xi(\tau) & = & 0 ; \\
\end{array}
\\
& \text{therefore, } [\xi,w(t)\partial_\tau]  =  B_1(t)w'(t)\partial_\tau.
\end{split}
\end{equation}
\\
The contribution from $v(t)\partial_t$ is given by,
\begin{equation}
\begin{split}
& \begin{array}{lcl}
\xi \: v(t)\partial_t(t) & = & \xi(v(t))=B_1(t)v'(t) \\
(-1)^{|\xi|}v(t)\partial_t \: \xi(t)& = & v(t)\partial_t(B_1(t))-v(t)\partial_t(B_2(t)) \\
& = & B_1(t)\tfrac{v(t)}{t} ; \\
\\
\xi \: v(t)\partial_t(\tau) & = & 0 \\
(-1)^{|\xi|}v(t)\partial_t \: \xi(\tau) & = & v(t)\partial_t(A_1(t))-v(t)\partial_t(A_2(t)) \\
& = & A_1(t)\tfrac{v(t)}{t} ; \\
\end{array}
\\
& \begin{array}{lcl}
\text{therefore, } [\xi,v(t)\partial_t] & = & \{B_1(t)v'(t)-B_1(t)\tfrac{v(t)}{t}\}\partial_t \\
& & -A_1(t)\tfrac{v(t)}{t}\partial_\tau.
\end{array}
\end{split}
\end{equation}
\\
The contribution from $m_0:=ad\tau-\tau^2\partial_\tau$ is given by,
\begin{equation} \label{eqn_mOform}
\begin{split}
& \begin{array}{lcl}
\xi m_0(t) & = & \xi(\tau t + t\tau) \\
& = & A(t)t + \tau(B_1(t)-B_2(t)) + B(t)\tau + t(A_1(t)-A_2(t)) \\
& = & 2tA_1(t) +[\tau,B_1(t)-B_2(t)] \\
(-1)^{|\xi|}m_0\xi(t) & = & m_0(B_1(t)-B_2(t))=[\tau,B_1(t)-B_2(t)] ; \\
\\
\xi m_0(\tau) & = & \xi(\tau^2)=A(t)\tau + \tau(A_1(t)-A_2(t)) \\
& = & [\tau,A_1(t)-A_2(t)] \\
(-1)^{|\xi|}m_0\xi(\tau)& = & m_0(A_1(t)-A_2(t))=[\tau,A_1(t)-A_2(t)] ; \\
\end{array}
\\
& \text{therefore, } [\xi,m_0] = 2tA_1(t)\partial_t.
\end{split}
\end{equation}

Adding the above equations together gives us the stated formula for $[\xi,m]$.
\end{proof}

If we assume that the element $2$ is invertible in the ground ring $R$ then by Theorem \ref{oddclass}, every odd Moore algebra is equivalent to one of the form $m_0+w(t)\partial_\tau$. This means that we need only calculate the Hochschild cohomology for odd Moore algebras of this form. Let us define $\widetilde{w}(t):=w(\sqrt{t})$.

\begin{prop} \label{nzhom}
Suppose we have an odd Moore algebra $A$ of the form $m_0+w(t)\partial_\tau$ where the element $2 \in R$ is invertible. If the first nonzero term of $w(t)$, which we denote by $w_k$, is not a zero divisor and $R$ has no $k$-torsion, then there is an isomorphism of $R$-modules $HH^*(A)\cong R[[t]]/(\widetilde{w}'(t))$.
\end{prop}

\begin{proof}
Lemma \ref{Hochformula} gives us the formula,
\begin{equation}
[\xi,m]=B_1(t)w'(t)\partial_\tau +2tA_1(t)\partial_t
\end{equation}
As $w_k$ is not a divisor of zero and $R$ has no $k$-torsion, $w'(t)=2kw_kt^{2k-1}+\ldots$ is not a divisor of zero. The Hochschild cocycles are then all the elements of the form $B(t)\partial_\tau + A(t)\partial_t$ where $A(t)$ and $B(t)$ are even power series with vanishing constant terms. The Hochschild coboundaries are all the elements of the form $B(t)\tfrac{w'(t)}{t} \partial_\tau + A(t)\partial_t$ where $A(t)$ and $B(t)$ are even power series with vanishing constant terms.

If $R_2[[t]]$ is the ring (without identity) of even power series with a vanishing constant term then the Hochschild cohomology is,
\[ HH^*(A)=R_2[[t]]/\dfrac{w'(t)}{t}R_2[[t]] \]
however, by the change of variable $t \mapsto t^2$ and the identity $w'(t)=2t\widetilde{w}'(t^2)$ we can obtain the isomorphism of $R$-modules $HH^*(A) \cong R[[t]]/(\widetilde{w}'(t))$.
\end{proof}

\begin{rem}
Since all the Hochschild cocycles are cohomologous to ones of the form $B(t)\partial_\tau$, the Lie bracket on $HH^*(A)$ is zero.
\end{rem}

Now I would like to turn my attention to the Hochschild cohomology of the trivial odd Moore algebra $R[X]/X^2$ which corresponds to choosing the characteristic power series $v(t)$ and $w(t)$ to be $0$. The \Ai-structure is then $m_0:=ad\tau-\tau^2\partial_\tau$. The following is a slightly different presentation of a standard result in homological algebra (cf. \cite[\S 9]{Weibel}).

\begin{prop} \label{zhom}
Assuming $2\in R$ is invertible, the Hochschild cohomology of the trivial odd Moore algebra $R[X]/X^2$ is the semi-direct sum of the Lie algebra $\{B(t)\partial_t, B(t) \ odd\}$ and the abelian Lie algebra $\{A(t)\partial_\tau, A(t) \ even\}$ such that,
\begin{displaymath}
[B(t)\partial_t,A(t)\partial_\tau]=B(t)A'(t)\partial_\tau
\end{displaymath}
where $A(t)$ and $B(t)$ have vanishing constant terms.
\end{prop}

\begin{proof}
Recalling equation \eqref{eqn_mOform} in Lemma \ref{Hochformula} we have,
\[ [\xi,m_0]=2tA_1(t)\partial_t \]
This means that the Hochschild cocycles are all the elements of the form $A(t)\partial_\tau + B(t)\partial_t$ where $A(t)$ is even and both $A(t)$ and $B(t)$ have vanishing constant terms. The Hochschild coboundaries are then all the elements of the form $B(t)\partial_t$ where $B(t)$ is even and has a vanishing constant term. The Hochschild cohomology is then given by,
\[ HH^*(R[X]/X^2)=\{A(t)\partial_\tau+B(t)\partial_t\} \]
where $A(t)$ is even, $B(t)$ is odd and both have vanishing constant terms.

Now to finish the proof we need only calculate $[B(t)\partial_t,A(t)\partial_\tau]$.
\begin{equation}
\begin{split}
& \begin{array}{lcl}
B(t)\partial_t \, A(t)\partial_\tau(\tau) & = & B(t)A'(t) \\
A(t)\partial_\tau \, B(t)\partial_t(\tau) & = & 0 ; \\
\\ 
\lbrack B(t)\partial_t,A(t)\partial_\tau \rbrack (t) & = & 0 ; \\
\end{array}
\\
& \text{therefore, } [B(t)\partial_t,A(t)\partial_\tau] = B(t)A'(t)\partial_\tau.
\end{split}
\end{equation}
\end{proof}

\begin{rem}
According to the theory developed in \cite{Penkava}, the even dimensional components (in the classical grading) of the Hochschild cohomology correspond to infinitesimal deformations of the \Ai-algebra by an even parameter. In the case of the trivial odd Moore algebra these are of the form $A(t)\partial_\tau$ where $A(t)$ is even. Motivated by the work done in \cite{Fuchs} and \cite{Alicemini}  it is reasonable to expect that these should provide a way to construct a ``miniversal'' deformation of the trivial Moore algebra and that, crudely speaking, the ``universal'' Moore algebra with $v(t)=0$ should be this miniversal deformation. This will be done in section \ref{Mini}.
\end{rem}
\section{Deformation of \Ai-Algebras} \label{Deform}

The purpose of this section is to set up the deformation theory for \Ai-algebras in a fairly general setting. We will assume that all \Ai-algebras are of the type described in definition \ref{def_AI}, that is to say that they are unital \Ai-algebras defined over an evenly graded ring $R$, in which the element $1 \in A$ can be completed to an $R$-basis. The \Ai-structure is then of the form \eqref{eqn_AI}.

The main differences between \cite{Penkava}, \cite{Alicemini} and the description that we will give here is that we will be using the cobar construction as opposed to the bar construction ($T\Sigma A$) in order to describe the deformation theory. In addition we will be considering the slightly more general situation of deforming an \Ai-algebra defined over a ring rather than just a field.

\begin{defi} \label{def_deform}
Suppose we have an \Ai-algebra $A$ with \Ai-structure $m$. Let $\Lambda$ be an (evenly graded) unital commutative $R$-algebra with an augmentation $\varepsilon:\Lambda \to R$. This $\varepsilon$ has a natural extension to a map of $\Lambda$-algebras $\varepsilon:\Lambda[[\tau,\boldsymbol{t}]] \to R[[\tau,\boldsymbol{t}]]$. A deformation of $A$ with base $(\Lambda,\varepsilon)$ is a unital, $\Lambda$-linear \Ai-structure $\bar{m}$ on the $\Lambda$-module $\Lambda \otimes_R A$ which renders the following diagram commutative:
\begin{displaymath}
\begin{CD}
\Lambda[[\tau,\boldsymbol{t}]] @>{\bar{m}}>> \Lambda[[\tau,\boldsymbol{t}]] \\
@V{\varepsilon}VV                                    @V{\varepsilon}VV \\
R[[\tau,\boldsymbol{t}]]       @>m>>         R[[\tau,\boldsymbol{t}]] \\
\end{CD}
\end{displaymath}
\end{defi}

\begin{defi} \label{def_equiv}
Suppose we have two deformations $\bar{m}_1$,$\bar{m}_2$ of an \Ai-algebra $A$ over the same base $(\Lambda,\varepsilon)$. An equivalence of deformations between $\bar{m}_1$ and $\bar{m}_2$ is a continuous automorphism $\phi$ of the $\Lambda$-algebra $\Lambda[[\tau,\boldsymbol{t}]]$ which renders the following diagrams commutative:
\begin{displaymath}
\xymatrix{\Lambda[[\tau,\boldsymbol{t}]] \ar^{\phi}[r] \ar_{\bar{m}_1}[d] &\Lambda[[\tau,\boldsymbol{t}]] \ar_{\bar{m}_2}[d] \\ \Lambda[[\tau,\boldsymbol{t}]] \ar^{\phi}[r] &\Lambda[[\tau,\boldsymbol{t}]]}
\qquad
\xymatrix{\Lambda[[\tau,\boldsymbol{t}]] \ar^{\varepsilon}[rr] \ar_{\phi}[rd] &&R[[\tau,\boldsymbol{t}]] \\ &\Lambda[[\tau,\boldsymbol{t}]] \ar_{\varepsilon}[ru]}
\end{displaymath}
\end{defi}

\begin{rem}
I shall refer to the automorphisms $\phi$ of the $\Lambda$-algebra $\Lambda[[\tau,\boldsymbol{t}]]$ which satisfy $\varepsilon \phi = \varepsilon$ as \emph{pointed automorphisms}. These \emph{pointed automorphisms} form a group under composition and act by conjugation on the left of the set consisting of deformations of $A$ over the base $(\Lambda,\varepsilon)$. Orbits of this action are in one to one correspondence with equivalence classes of deformations.
\end{rem}

I shall now define the deformation functor for an \Ai-algebra. Given two augmented $R$-algebras $(\Lambda,\varepsilon)$ and $(\Lambda',\varepsilon')$, a homomorphism of augmented $R$-algebras $f:\Lambda \to \Lambda'$ is a (unital) homomorphism of $R$-algebras such that $\varepsilon'f=\varepsilon$.

\begin{defi} \label{def_push-out}
Suppose we have an \Ai-algebra $A$ defined over an evenly graded and commutative ring $R$. The deformation functor is defined on the category of unital augmented $R$-algebras by sending the algebra $(\Lambda,\varepsilon)$ to the set of equivalence classes of all possible deformations of $A$ with base $(\Lambda,\varepsilon)$.

To define it on morphisms suppose we have a homomorphism of augmented $R$-algebras $f:\Lambda \to \Lambda'$ and a deformation of $A$ with base $(\Lambda,\varepsilon)$, that is an \Ai-structure $\bar{m}_\Lambda$ as described in definition \ref{def_deform}. $f$ naturally extends to a homomorphism of $\Lambda$-algebras $f:\Lambda[[\tau,\boldsymbol{t}]] \to \Lambda'[[\tau,\boldsymbol{t}]]$. We must define a deformation $\bar{m}_{\Lambda'}$ of $A$ with base $(\Lambda',\varepsilon')$ which we will refer to as the push-out of $\bar{m}_\Lambda$ by $f$:
\begin{displaymath}
\bar{m}_{\Lambda'}:=(f\bar{m}_\Lambda(\tau))\partial_\tau +\sum_{i \in I}(f\bar{m}_\Lambda(t_i))\partial_{t_i}
\end{displaymath}
Clearly the push-out operation respects equivalences of deformations and hence gives rise to a map between equivalence classes of deformations.
\end{defi}

\begin{rem}
The deformation theory of \Ai-algebras is a generalisation of the theories developed in \cite{Murray} and \cite{dgadef}. In \cite{Murray} Gerstenhaber considered deformations of associative algebras and related infinitesimals to Hochschild cohomology classes. In describing this generalisation we will need to make use of the alternative definition of \Ai-algebras involving higher multiplication maps as described in \cite{Bernhard} and \cite{Andrey}. An associative graded algebra $m_2:A \otimes A \to A$ can be viewed as an \Ai-algebra by specifying all the other multiplications $m_i:A^{\otimes i} \to A$ ($i \neq 2$) to be zero. An isomorphism $f_1:A_1\to A_2$ of associative graded algebras can be viewed as an isomorphism of \Ai-algebras by specifying the maps $f_i:A_1^{\otimes i} \to A_2$ ($i \geq 2$) to be zero. Let $F$ be the deformation functor for \Ai-algebras defined above and let $G$ be the deformation functor for associative graded algebras which sends augmented $R$-algebras $\Lambda$ to equivalence classes of deformations of the associative graded algebra $A$ with base $\Lambda$. The above observations then give us an injective morphism of functors $G \to F$ which describes the way in which the deformation theory for \Ai-algebras encompasses the deformation theory for associative algebras. 

In a similar manner a DGA $m_1:A \to A$, $m_2:A\otimes A \to A$ can be regarded as an \Ai-algebra by specifying all higher multiplications $m_i:A^{\otimes i} \to A$ ($i\geq 3$) to be zero. An isomorphism $f:A_1 \to A_2$ of DGA's can be viewed as an isomorphism of \Ai-algebras by specifying the maps $f_i:A_1^{\otimes i} \to A_2$ ($i \geq 2$) to be zero. Now if $G$ is the deformation functor for DGA's defined by sending the augmented $R$-algebra $\Lambda$ to equivalence classes of deformations of the DGA $A$ with base $\Lambda$ and $F$ is the deformation functor for \Ai-algebras as above, then the preceding observations give us a morphism of functors $G \to F$.
\end{rem}

Now we must define the appropriate notions of versality, universality and miniversality for deformations. In deformation theory in general it is not always possible to find a universal deformation of a given object. It is, however, much more likely that a ``miniversal'' deformation exists. Indeed in \cite{Fuchs} it was shown that for Lie algebras there is, under certain restrictions, an explicit construction for this miniversal deformation in terms of the Lie algebra's cohomology. A similar method was used in \cite{Alicemini} to show that under certain restrictions, miniversal deformations exist for \Ai-algebras.

\begin{defi} \label{def_mini}
Suppose we have an \Ai-algebra $A$ and a deformation $\bar{m}_\Lambda$ of $A$ with base $(\Lambda,\varepsilon)$. $\bar{m}_\Lambda$ is a versal deformation of $A$ if given any other deformation $\bar{m}_{\Lambda'}$ of $A$ with base\footnote{In \cite{Fuchs} the restriction was made that these bases be local, however, as we are working over rings and not just fields, this restriction is unsuitable.} $(\Lambda',\varepsilon')$ there is a homomorphism $f:\Lambda \to \Lambda'$ such that $\bar{m}_{\Lambda'}=f_*(m_\Lambda)$. That is to say that, up to an equivalence of deformations, $\bar{m}_{\Lambda'}$ is the push-out of $\bar{m}_\Lambda$ by $f$.

If in addition the homomorphism $f$ is unique for any deformation of $A$, the deformation $\bar{m}_\Lambda$ is called \emph{universal}. The deformation $\bar{m}_{\Lambda'}$ is called \emph{infinitesimal} if $\ker(\varepsilon')^2=0$. If the homomorphism $f$ is unique for any \emph{infinitesimal} deformation of $A$, the deformation $\bar{m}_\Lambda$ is called \emph{miniversal} (cf. \cite{Fuchs}, \cite{Alicemini}).
\end{defi}
\section{The Miniversal Deformations of $R[X]/X^2$} \label{Mini}

In this section we will calculate the miniversal deformations of the trivial Moore algebra $R[X]/X^2$ and show that there is no universal deformation of this algebra. In \cite[\S 5]{Andrey} the universal even and odd Moore algebras were constructed. The universal even Moore algebra of degree $d$ is the Moore algebra defined over the ring $\mathbb{Z}[u_1,u_2,\ldots]$ where $|u_i|=i(d+2)-2$ which is given by the \Ai-structure,
\begin{equation} \label{eqn_evenuni}
m:=\sum_{i=1}^\infty u_it^i\partial_\tau + ad\tau-\tau^2\partial_\tau
\end{equation}
Let us define the $R$-universal even Moore algebra of degree $d$ to be the even Moore algebra defined over the ring $R[u_1,u_2,\ldots]$ with the same \Ai-structure as \eqref{eqn_evenuni}. This is the universal object in the sense of \cite[\S 5]{Andrey} only for even Moore algebras defined over a commutative $R$-algebra.

The odd universal Moore algebra of degree $d$ is the Moore algebra defined over the ring $\mathbb{Z}[v_1,v_2,\ldots]\otimes\mathbb{Z}[w_1,w_2,\ldots]$ where $|v_i|=2i(d+2)-d-3$ and $|w_i|=2i(d+2)-2$ which is given by the \Ai-structure \eqref{eqn_odddef}. Let us define the $R$-universal odd Moore algebra of degree $d$ with $v(t)=0$ to be the odd Moore algebra defined over the ring $R[w_1,w_2,\ldots]$ with \Ai-structure,
\begin{equation} \label{eqn_odduni}
m:=\sum_{i=1}^\infty w_it^{2i}\partial_\tau + ad\tau -\tau^2\partial_\tau
\end{equation}

\begin{theorem} \label{evenmini}
The $R$-universal even Moore algebra is a miniversal deformation of the trivial even Moore algebra $A:=R[X]/X^2$.
\end{theorem}

\begin{proof}
Clearly the $R$-universal even Moore algebra is a deformation of $A$ with base $R[u_1,u_2,\ldots]$. Suppose we have an arbitrary deformation $\bar{m}$ of $A$ with base $(\Lambda,\varepsilon)$. The \Ai-structure $\bar{m}$ must have the form,
\[ \bar{m}=\sum_{i=1}^\infty \lambda_it^i\partial_\tau + m_A \]
where $\lambda_i \in \Lambda$, $|\lambda_i|=i(d+2)-2$ and $m_A:=ad\tau-\tau^2\partial_\tau$ is the (extension of the) \Ai-structure on $A$. The $R$-algebra $\Lambda$ has a decomposition $\Lambda=R\oplus\ker\varepsilon$ and the restrictions of definition \ref{def_deform} mean that $\lambda_i \in \ker\varepsilon$. Choose a homomorphism $f:R[u_1,u_2,\ldots] \to \Lambda$ by sending $u_i \mapsto \lambda_i$. The push-out of \eqref{eqn_evenuni} by $f$ is $\bar{m}$ so the $R$-universal even Moore algebra is a versal deformation of $A$.

Let us show that this is actually a miniversal deformation of $A$. For this purpose let us assume that $\ker(\varepsilon)^2=0$. An equivalence of deformations for deformations over the base $(\Lambda,\varepsilon)$ is given by a pointed automorphism of the $\Lambda$-algebra $\Lambda[[\tau,t]]$ corresponding to power series $G(t),F(t) \in \Lambda[[t]]$ with vanishing constant terms where $G(t)=0$ and $F(t)=f_1t^1+f_2t^2+\ldots$ ($f_1$ invertible). Since this automorphism is pointed we must have $\varepsilon(f_1)=1$ and $\varepsilon(f_i)=0$ for $i\ge2$. These pointed automorphisms act on $\bar{m}$ by $\lambda(t) \mapsto \lambda(F(t))$ (cf. \cite[\S 6]{Andrey}), however since $\lambda_i \in \ker\varepsilon$, $\varepsilon(f_1)=1$, $f_i \in \ker\varepsilon$ for $i\ge2$ and $\ker(\varepsilon)^2=0$, this is a trivial action. It therefore follows that the map $f$ must be unique.
\end{proof}

\begin{theorem} \label{oddmini}
Assuming the element $2 \in R$ is invertible, the $R$-universal odd Moore algebra with $v(t)=0$ is a miniversal deformation of the trivial odd Moore algebra $A:=R[X]/X^2$.
\end{theorem}

\begin{proof}
Clearly the $R$-universal odd Moore algebra with $v(t)=0$ is a deformation of $A$ with base $R[w_1,w_2,\ldots]$. Suppose we have an arbitrary deformation $\bar{m}$ of $A$ with base $(\Lambda,\varepsilon)$. The \Ai-structure must have the form,
\[ \bar{m}=\sum_{i=1}^\infty v_it^{2i}\partial_t +\sum_{i=1}^\infty u_it^{2i}\partial_\tau +m_A \]
In addition we must have $v_i,u_i \in \ker \varepsilon$. This means we can choose a pointed automorphism $\phi:=(-\tfrac{v(t)}{2t},t)$ which gives us by \eqref{eqn_oddact} an equivalent deformation of the form $\sum_{i=1}^\infty \lambda_it^{2i}\partial_\tau+m_A$ where $\lambda_i \in \ker\varepsilon$. Choose a homomorphism $f:R[w_1,w_2,\ldots] \to \Lambda$ by sending $w_i \mapsto \lambda_i$. The push-out of \eqref{eqn_odduni} by $f$ is equivalent to $\bar{m}$ so the $R$-universal odd Moore algebra with $v(t)=0$ is a versal deformation of $A$.

To prove that this is a miniversal deformation, suppose that $\bar{m}$ is an infinitesimal deformation of $A$. We need to show that if two deformations,
\begin{displaymath}
\sum_{i=1}^\infty \lambda_it^{2i}\partial_\tau +m_A \quad , \quad \sum_{i=1}^\infty \lambda'_it^{2i}\partial_\tau +m_A
\end{displaymath}
are equivalent ($\lambda_i , {\lambda'}_i \in \ker \varepsilon$) then they are equal.

Suppose they are equivalent by an equivalence of deformations  $\phi=(G(t),F(t))$, then \eqref{eqn_oddact} implies that $G(t)=0$. $F(t)$ then acts according to \eqref{eqn_oddact} by sending $\lambda(t) \mapsto \lambda(F(t))$, however, we have already shown in Theorem \ref{evenmini} that this is a trivial action and therefore the two deformations are equal.
\end{proof}

Finally we would like to know whether or not it is possible to construct universal deformations of the trivial Moore algebra $R[X]/X^2$. This turns out to be impossible.

\begin{prop} \label{nonuniv}
There do not exist universal deformations of either the odd or the even trivial Moore algebra $A:=R[X]/X^2$.
\end{prop}

\begin{proof}
We shall treat only the even case, the odd case being virtually identical. First of all let us consider a deformation $\bar{m}$ of $A$ with base $R[x,y]$ ($|x|:=d$ and $|y|:=0$) given by $\bar{m}:=xt\partial_\tau+m_A$. Using the pointed automorphism $\phi:=(0,t+yt)$ we see that this deformation is equivalent to $(x+xy)t\partial_\tau+m_A$. These clearly come from two different push-outs $f_1$ and $f_2$ of the miniversal deformation of $A$.

Now we must use the functoriality of the push-out operation and the fact that it respects equivalences of deformations. If a universal deformation $m_U=\sum_{i=1}^\infty \lambda_it^i\partial_\tau + m_A$ of $A$ over a base $\Lambda$ exists, then there is a homomorphism $g:\Lambda \to R[u_1,u_2,\ldots]$ such that \eqref{eqn_evenuni} is (up to equivalence) the push-out of $m_U$ by $g$. Universality of $m_U$ then implies that $f_1g=f_2g$.

Consider the infinitesimal deformation of $A$ over the base $R[u_1,u_2,\ldots]/(u_iu_j)$ which is given by the derivation $\sum_{i=1}^\infty u_it^i\partial_\tau + m_A$. This is the push-out of \eqref{eqn_evenuni} by the natural quotient map $\sigma:R[u_1,u_2,\ldots] \to R[u_1,u_2,\ldots]/(u_iu_j)$ and is therefore equivalent to the push-out of $m_U$ by $\sigma g$, however, in Theorem \ref{evenmini} it was shown that pointed automorphisms act trivially on infinitesimal deformations, therefore
\[ g(\lambda_k)=u_k\mod(u_iu_j), \qquad k \geq 1\]
It is clear from the construction of $f_1$ and $f_2$ and the equation $f_1g(\lambda_1)=f_2g(\lambda_1)$ that we have a contradiction.
\end{proof}
\section{Formal One-Parameter Deformations of \Ai-Algebras} \label{Oneparam}

This subject has already been touched upon by Penkava and Weldon in \cite{Weldon}. Formal one-parameter deformations of DGA's were considered in \cite{dgadef}, however there are certain advantages to working with \Ai-algebras as opposed to ordinary DGA's. In particular it was only possible in \cite{dgadef} to show the link between infinitesimals of deformations of a DGA and the corresponding cohomology theory for DGA's whose ground ring contained $\mathbb{Q}$. With \Ai-algebras however, there are no such problems.

The following work draws a parallel between the deformation theories described in \cite{Weldon} and \cite{dgadef}. The theory of formal one-parameter deformations of an \Ai-algebra $A$ expounded here is parallel to the theory of formal one-parameter deformations of the DGA $R[[\tau,\boldsymbol{t}]]$ which describes the \Ai-structure on $A$. We will therefore need to make use of the work done by Gerstenhaber and Wilkerson on the subject in \cite{dgadef}.

The ring of formal power series $R[[s]]$, where $s$ has even degree, has a fundamental system of neighbourhoods which are powers of the augmentation ideal $R[[s]]/R$. This fundamental system gives a topology on the $R[[s]]$-module $A[[s]]$. A formal one-parameter deformation $m_s$ of $A$ is given by a continuous, $R[[s]]$-linear \Ai-structure on the $R[[s]]$-module $A_s:=A[[s]]$.

\begin{defi}
Suppose we have an \Ai-algebra $A$ defined over an evenly graded and commutative ring $R$ with \Ai-structure $m$. A formal one-parameter deformation $m_s$ of $A$ is given by a continuous (both in the sense of definition \ref{def_AI} and in the sense of the topology inherited from $R[[s]]$) $R[[s]]$-linear derivation on the $R[[s]]$-algebra $R[[s]][[\tau,\boldsymbol{t}]]$. Such a derivation $m_s$ is any derivation of the form,
\[ m_s=m+sm_1+s^2m_2+\ldots+s^nm_n+\ldots \]
where $m_i$ is (the extension of) a normalised derivation from $R[[\tau,\boldsymbol{t}]]$ to $R[[\tau,\boldsymbol{t}]]$ which has degree $i|s|+2$ (in the classical grading). By convention we define $m_0:=m$. Furthermore we require this map to give an \Ai-structure on $A[[s]]$, that is to say that $m_s^2=0$ or equivalently,
\begin{equation} \label{eqn_defcri}
\sum_{\begin{subarray}{c} i+j=k \\ i,j > 0 \end{subarray}}m_im_j=-[m_k,m] , \quad \text{for all } k \geq 1
\end{equation}
\end{defi}

\begin{defi}
A formal one-parameter automorphism $\phi_s$ (the analogue of a pointed automorphism) is given by a continuous, $R[[s]]$-linear endomorphism of $R[[s]][[\tau,\boldsymbol{t}]]$. Such an endomorphism $\phi_s$ is any endomorphism of the form,
\[ \phi_s=1+s\phi_1+s^2\phi_2+\ldots+s^n\phi_n+\ldots \]
where $\phi_i$ is (the extension of) a continuous $R$-linear map from $R[[\tau,\boldsymbol{t}]]$ to $R[[\boldsymbol{t}]]/R$. Furthermore we require that $\phi_s(1)=1$ or equivalently that $\phi_i(1)=0$ for all $i\geq 1$. We also make the restriction that $\phi_s$ be multiplicative, a condition that can be expressed as,
\begin{equation} \label{eqn_autcri}
\sum_{\begin{subarray}{c} i+j=k \\ i,j \geq 0 \end{subarray}}\phi_i(a)\phi_j(b)=\phi_k(ab), \quad \text{for all } a,b \in R[[\tau,\boldsymbol{t}]]
\end{equation} 
\end{defi}

\begin{defi}
Let $m_s$ and $m'_s$ be two formal one-parameter deformations of an \Ai-algebra $A$. We say that $m_s$ is equivalent to $m'_s$ if there exists a formal automorphism $\phi_s$ such that $m'_s=\phi_sm_s\phi_s^{-1}$.
\end{defi}

\begin{rem}
The formal one-parameter automorphisms form a group under composition and act by conjugation on the left of the space of formal one-parameter deformations of $A$. Orbits of this action are in one to one correspondence with equivalence classes of formal one-parameter deformations of $A$.
\end{rem}

\begin{rem}
There is a similar definition of a formal one-parameter deformation of an \Ai-algebra $A$ in terms of higher multiplication maps. The restrictions on the maps $m_i:A[[s]]^{\otimes i} \to A[[s]]$ are the same as those alluded to in remark \ref{rem_hmult} except that the $m_i$'s must be continuous. This definition is entirely consistent with the definition given above however.
\end{rem}

Suppose we have a deformation $m_s=m+s^km_k+s^{k+1}m_{k+1}+\ldots$ , then \eqref{eqn_defcri} shows that $[m_k,m]=0$, that is to say that $m_k$ is a Hochschild cocycle. Similarly, given a formal automorphism $\phi_s=1+s^k\phi_k+s^{k+1}\phi_{k+1}+\ldots$ , \eqref{eqn_autcri} determines that $\phi_k$ must be a (normalised) derivation of $R[[\tau,\boldsymbol{t}]]$ (note that since $\phi_k$ is even, a derivation is the same thing as a graded derivation).

We say that $\phi_s=1+s\phi_1+\ldots+s^n\phi_n$ is a formal automorphism of order $n$ if \eqref{eqn_autcri} holds for all $1\leq k\leq n$. Similarly $m_s=m+sm_1+\ldots+s^nm_n$ is a deformation of order $n$ if \eqref{eqn_defcri} holds for all $1\leq k\leq n$. In order to show the link between Hochschild cohomology and infinitesimals of deformations we will need the following relatively trivial lemma:

\begin{lemma} \label{extlem}
Suppose we are given a formal automorphism of order $n$,
\[ \psi_s=1+s\phi_1+\ldots+s^n\phi_n \]
then we can extend $\psi_s$ to a formal automorphism,
\[ \phi_s=1+s\phi_1+\ldots+s^n\phi_n+s^{n+1}\phi_{n+1}+s^{n+2}\phi_{n+2}+\ldots\]
by specifying the action of $\phi_{n+1},\phi_{n+2},\ldots$ on the generators $\tau,t_1,t_2,\ldots \ $. Furthermore such a choice uniquely determines $\phi_s$.
\end{lemma}

\begin{proof}
Suppose that we are given a choice $\phi_i(\tau),\phi_i(t_1),\phi_i(t_2),\ldots$ for the action of $\phi_{n+1},\phi_{n+2},\ldots$ on the generators (note that if the collection $\boldsymbol{t}=t_1,t_2,\ldots$ is infinite then this will impose some restrictions on the collection $\phi_i(\tau),\phi_i(t_1),\phi_i(t_2),\ldots$ ). We can extend $\phi_{n+1},\phi_{n+2},\ldots$ to continuous $R$-linear maps from $R[[\tau,\boldsymbol{t}]]$ to $R[[\boldsymbol{t}]]$ in the following way,
\[ \phi_k(x_1\ldots x_m)=\sum_{i_1+\ldots+i_m=k}\phi_{i_1}(x_1)\ldots\phi_{i_m}(x_m) \]
where $x_1,\ldots,x_m \in \{\tau,t_1,t_2,\ldots\}$. It is then straightforward to check that the $\phi_k$'s so defined satisfy \eqref{eqn_autcri}. Clearly the map $\phi_s$ is uniquely determined by its action on the generators which is in turn determined by the action of $\phi_{n+1},\phi_{n+2},\ldots$ on the generators.
\end{proof}

One of the problems with the deformation theory for an ordinary DGA is that this extension may not always be possible. If it is however, then Gerstenhaber and Wilkerson showed in \cite{dgadef} that infinitesimals of deformations correspond to odd cohomology classes (the cohomology theory used by Gerstenhaber and Wilkerson was not Hochschild cohomology, our infinitesimals will live in even components of the classically graded Hochschild cohomology).

\begin{defi}
Suppose we have an \Ai-algebra $A$ with \Ai-structure $m$. An infinitesimal automorphism of $A$ is a formal automorphism $\phi_s:=1+s^k\phi_k$ of order $k$ where $\phi_k$ is a Hochschild cocycle of degree $k|s|+1$ (in the classical grading). We say that $\phi_s$ is integrable if it can be extended to a formal automorphism $\psi_s=1+s^k\phi_k+s^{k+1}\phi_{k+1}+\ldots$ such that $[\psi_s,m]=0$, that is to say that $\psi_s$ is a continuous automorphism of the \Ai-algebra $A[[s]]$ whose \Ai-structure is given by the deformation $m_s:=m$.
\end{defi}

\begin{rem}
Note that if $\mathbb{Q} \subset R$ then every infinitesimal automorphism is integrable. Since we can exponentiate we can choose the extension of $\phi_s$ to be $\psi_s:=\exp(s^k\phi_k)$. The problem of integrating an infinitesimal automorphism of an associative algebra when working over a field of characteristic $p$ was known to Gerstenhaber and was discussed in \cite{Murray}. There are known to be examples of infinitesimal automorphisms of an associative algebra which cannot be integrated when working over a field of nonzero characteristic.
\end{rem}

The method of proof of the following theorem borrows from the work done by Gerstenhaber and Wilkerson in \cite{dgadef}. Let us call a deformation \emph{trivial} if it is equivalent to the deformation $m_s:=m$.

\begin{theorem} \label{Hochinf}
\
\begin{enumerate}
\item[(i)]
Suppose we have an \Ai-algebra $A$ with \Ai-structure $m$. Every formal one-parameter deformation of $A$ is either trivial or equivalent to a deformation of the form,
\[ m_s=m+s^km_k+s^{k+1}m_{k+1}+\ldots \]
where $m_k$ is a Hochschild cocycle of degree $k|s|+2$ (in the classical grading) which is not cohomologous to zero. We will refer to $m_k$ as the infinitesimal of the deformation.

\item[(ii)]
Suppose that every infinitesimal automorphism of $A$ is integrable and that we have two equivalent deformations,
\[ m_s=m+s^km_k+\ldots \quad , \quad m'_s=m+s^nm'_n+\ldots \qquad (n\geq k)\]
where $m_k$ is not cohomologous to zero, then $n=k$ and $m_k$ is cohomologous to $m'_n$.
\end{enumerate}
\end{theorem}

\begin{rem}
One should not confuse the notion of an infinitesimal of a deformation with that of infinitesimal deformations. Remember an infinitesimal deformation is a deformation over a base whose augmentation ideal has square equal to zero.
\end{rem}

\begin{proof}
\
\begin{enumerate}
\item[(i)]
First of all let us show that if we have a deformation $m_s=m+s^km_k+\ldots$ then we can transform $m_k$ to any cohomologous cocycle using a formal automorphism. Suppose we are given a normalised derivation $\xi$ (of classical degree $k|s|+1$), then by Lemma \ref{extlem} we can extend $1+s^k\xi$ to a formal automorphism,
\[ \phi_s=1+s^k\xi+s^{k+1}\phi_{k+1}+\ldots \]
As $\phi_s^{-1}=1-s^k\xi+\ldots$ we have $\phi_sm_s\phi_s^{-1}=m+s^k(m_k+[\xi,m])+\ldots \ $.

\paragraph*{}
Now suppose we have a deformation $m_s=m+sm_1+s^2m_2+\ldots$ such that all equivalent deformations begin with a cocycle cohomologous to zero. Using the above procedure we can find a map $\phi_1:R[[s,\tau,\boldsymbol{t}]] \to R[[s,\boldsymbol{t}]]$ such that $(1+s\phi_1)m_s(1+s\phi_1)^{-1}$ is a deformation $m'_s$ of the form $m+s^2m'_2+\ldots \ $. Similarly we can find a map $\phi_2$ such that $(1+s^2\phi_2)m'_s(1+s^2\phi_2)^{-1}$ has the form $m+s^3m_3+\ldots$ and construct a sequence of maps $\phi_1,\phi_2,\phi_3\ldots$ such that
\[ m=\ldots(1+s^2\phi_2)(1+s\phi_1)m_s(1+s\phi_1)^{-1}(1+s^2\phi_2)^{-1}\ldots \]
This means that the deformation $m_s$ must be trivial.

\item[(ii)]
Suppose that every infinitesimal automorphism of $A$ is integrable and that we have two equivalent deformations,
\[ m_s=m+s^km_k+\ldots \quad , \quad m'_s=m+s^nm'_n+\ldots \qquad (n \geq k) \]
where $m_k$ is not cohomologous to zero. There is a formal automorphism $\phi_s$ such that $\phi_sm'_s\phi_s^{-1}=m_s$. We will show that we can write $\phi_s$ in the form,
\begin{equation} \label{eqn_formalform}
\phi_s=(1+s^k\phi_k+\ldots)\psi_s[k-1]\ldots\psi_s[1]
\end{equation}
where $\psi_s[i]=1+s^i\psi_i+\ldots$ is a formal automorphism which commutes with $m$. First of all if $k=1$ then \eqref{eqn_formalform} is vacuous so we can assume $k>1$. Writing $\phi_s$ as a series $\phi_s=1+s\phi_1+\ldots$ we see that since $\phi_sm'_s\phi_s^{-1}=m_s$, we have $[\phi_1,m]=0$. Next integrate the infinitesimal automorphism $1+s\phi_1$ to a formal automorphism $\psi_s[1]$ commuting with $m$ and express $\phi_s$ as $\phi_s=(1+s^2\phi_2+\ldots)\psi_s[1]$. If $k>2$ then $[\phi_2,m]=0$ and we can proceed to construct $\psi_s[i]$ inductively. Now that we have $\phi_s$ in the form of \eqref{eqn_formalform} and since $\phi_sm'_s\phi_s^{-1}=m_s$ we can conclude that $m_k=m'_k+[\phi_k,m]$, hence $k=n$ and $m_k$ is cohomologous to $m'_n$.
\end{enumerate}
\end{proof}

\begin{rem}
Part (i) of Theorem \ref{Hochinf} is the analogue of well known results in algebraic deformation theory (cf. \cite{Murray}, \cite{dgadef}), however, part (ii) contains a result on the deformation theory of associative algebras which may be new.

There are known to be examples of associative algebras in which there exist infinitesimal automorphisms which cannot be integrated and for which part (ii) of Theorem \ref{Hochinf} does not hold (cf. \cite[\S 6]{Murray2}). That is to say that there are trivial deformations which begin with a Hochschild cocycle noncohomologous to zero. Perhaps the simplest example is provided by the associative algebra $A:=(\mathbb{Z}/2\mathbb{Z})[X]/(X^2)$. Choose a formal automorphism $\phi_s$ by sending $X$ to $X+s\cdot 1$ , then $\phi_s$ acts on the deformation $m'_s:=m$ to give the following deformation:
\[ m_s(x,x)=s^2 \cdot 1 \]
Since we are working in characteristic 2, the differential $\delta$ is zero on normalised Hochschild cochains, therefore although $m_s$ is trivial, it begins with a Hochschild cocycle noncohomologous to zero.
\end{rem}

We say that an \Ai-algebra $A$ with \Ai-structure $m$ is \emph{rigid} if every formal one-parameter deformation of it is trivial. In view of part (i) of Theorem \ref{Hochinf} we have the following corollary:

\begin{cor} \label{cor_rigid}
Suppose we have an \Ai-algebra $A$ with $HH^{2*}(A)=0$ (here we use the classical grading), then $A$ is rigid.
\end{cor}

\begin{rem}
One advantage of working with \Ai-algebras is that although Gerstenhaber and Wilkerson were able to prove the analogue of Corollary \ref{cor_rigid} for DGA's whose ground ring contained $\mathbb{Q}$ (cf. \cite{dgadef}), this restriction is not necessary when deforming \Ai-algebras. 
\end{rem}

Finally we would like to know how to extend a deformation $m_s=m+sm_1+\ldots+s^nm_n$ of order $n$ to a deformation of order $n+1$. The obstruction to this extension $\Obs(m_s):=\sum_{i+j=n+1}m_im_j$ ($i,j>0$) turns out to be an odd Hochschild cohomology class (in the classical grading). The method of proof of the following theorem was used by Gerstenhaber and Wilkerson in \cite{dgadef} in which they investigated formal one-parameter deformations of differential graded modules.

\begin{theorem}
Suppose we have a deformation $m_s=m+sm_1+\ldots+s^nm_n$ of order $n$ of an \Ai-algebra $A$, then $\Obs(m_s)$ is an odd Hochschild cocycle (in the classical grading). The deformation $m_s$ is extendible to a deformation of order $n+1$ if and only if $\Obs(m_s)$ is cohomologous to zero.
\end{theorem}

\begin{proof}
First of all let us show that $\Obs(m_s)$ is a normalised derivation. We have,
\[ m_im_j(ab)=m_im_j(a)b+(-1)^{|a|}(m_i(a)m_j(b)-m_j(a)m_i(b))+am_im_j(b) \]
however the middle term sums to zero in $\sum_{i+j=n+1}$ and hence $\Obs(m_s)$ is a graded derivation (in the standard grading).

Next let us show that $\Obs(m_s)$ is a Hochschild cocycle. Since \eqref{eqn_defcri} holds for all $k\leq n$ we have,
\[ [m_k,m]=-\sum_{i+j=k}m_im_j \quad (i,j>0) \]
As the map $ad \, m:\xi \mapsto [m,\xi]$ is a graded derivation (in the standard grading) on the algebra of $R$-linear maps $\End_R(R[[\tau,\boldsymbol{t}]])$ we have,
\begin{displaymath}
\begin{array}{lcl}
[m,\Obs(m_s)] & = & \sum_{i+j=n+1}[m,m_i]m_j-\sum_{i+j=n+1}m_i[m,m_j] \\
 & = & -\sum_{i+j+k=n+1}m_im_jm_k+\sum_{i+j+k=n+1}m_im_jm_k = 0 \\
\end{array}
\end{displaymath}
where $i,j,k>0$. We can find a normalised derivation $m_{n+1}$ such that \eqref{eqn_defcri} holds for $k=n+1$ precisely when $\Obs(m_s)$ is a Hochschild coboundary.
\end{proof}

In view of the following proposition it would seem more appropriate to interpret the obstruction as a cohomology class.

\begin{prop} \label{prop_obsequiv}
If two deformations $m_s$ and $m'_s$ of order $n$ are equivalent, then they have cohomologous obstructions.
\end{prop}

\begin{proof}
First of all notice that deformations and formal automorphisms of order $n$ are really deformations and pointed automorphisms over the base $R[s]/(s^{n+1})$. Suppose we have a formal automorphism $\psi_s$ of order $n$ such that $m'_s=\psi_sm_s\psi_s^{-1}$ and extend it to a formal automorphism $\phi_s$ using Lemma \ref{extlem}.

Considering $m_s$ as an $R[[s]]$-linear derivation on $R[[s]][[\tau,\boldsymbol{t}]]$ we have,
\[ m_s^2=s^{n+1}\Obs(m_s)+O(s^{n+2}) \]
Since $\phi_s \equiv \psi_s \mod (s^{n+1})$ we have $m'_s=\phi_sm_s\phi_s^{-1}+s^{n+1}\xi+O(s^{n+2})$, where $\xi$ is a normalised derivation of even degree (in the classical grading). This means that,
\begin{displaymath}
\begin{array}{lcl}
{m'_s}^2 & = & \phi_sm_s^2\phi_s^{-1}+s^{n+1}[\xi,m]+O(s^{n+2}) \\
 & = & s^{n+1}(\Obs(m_s)+[\xi,m])+O(s^{n+2}) \\
\end{array}
\end{displaymath}
and therefore $\Obs(m'_s)=\Obs(m_s)+[\xi,m]$.
\end{proof}

\begin{rem}
Notice that in order to prove this result it was necessary to use Lemma \ref{extlem}. This suggests that whilst true for \Ai-algebras, this result might not hold for DGA's in general. Indeed in \cite{dgadef} Gerstenhaber and Wilkerson state this result for differential graded modules, but make no mention as to its truth for DGA's.
\end{rem}

\end{document}